\font\twbf=cmbx12
\magnification=1200     \hfuzz=1pt
\hsize=125mm \vsize=187mm \hoffset=4mm \voffset=10mm
\parindent=0mm
\def\parag{\vk{2}}

\def\large{\let\rm=\xii \let\bf=\twbf \let\sl=\twsl \let\sf=\twsf
\let\sfi=\twsfi \baselineskip=14pt minus1pt \rm}
\def\normal{\let\rm=\tenrm \let\bf=\tenbf \let\sl=\tensl \let\sf=\tensf
\let\sfi=\tensfi \baselineskip=12pt minus1pt \rm}
\def\small{\let\rm=\viii \let\bf=\eightbf \let\sl=\eightsl \let\sf=\eightsf
\let\sfi=\eightsfi \baselineskip=9.5pt minus.75pt \rm}
\def\tiny{\let\rm=\vii \let\bf=\sevenbf \let\sl=\sevensl \let\sf=\sevensf
\let\sfi=\sevensfi \baselineskip=8.3pt minus.65pt \rm}

\def\fdate#1#2#3{{\oldstyle #1{\rm~#2~}#3}}

\def\date{{\oldstyle\the\day}\ \ifcase\month\or janvier\or f\'evrier\or mars
\or avril\or mai\or juin\or juillet\or ao\^ut\or septembre\or octobre
\or novembre\or d\'ecembre\fi\ {\oldstyle\the\year}}

\def\qed{\hfill\kern 6pt
\lower 2pt\hbox{\vrule\vbox to10pt{\hrule width 4pt \vfil\hrule}\vrule}\par}

\def\cf{{\sl cf.\/}\ }

\parskip=0pt plus 1pt minus 1pt

\def\und#1{\ifmmode\underline{#1}\else\underbar{#1}\fi}

\def\,{\mkern2mu}
\def\;{\mkern4mu}

\def\im{{\rm i}}
\def\op#1{\mathop{\rm #1}}
\def\Log{\op{Log}}
\def\exp{\op{exp}}
\def\e#1{\op{e}\nolimits^{#1}}
\def\ei#1{\op{e}\nolimits^{{\im}#1}}

\def\Lim{\op{Lim}}
\def\Sup{\op{Sup}}

\def\dif{\not=}
\def\0{\hb{\rm\O}}
\def\1{[0,1]}

\def\cC{{\cal C}}
\def\cD{{\cal D}}

\def\R{{\rm I\! R}}
\def\E{{\rm I\! E}}
\def\Frac#1#2{\disp{{#1}\over {#2}}}
\def\rond{\ifmmode{\mathrel{\circ}}\else{$\circ$}\fi}

\def\parag{\par\ \par}

\def\hb{\hbox}
\def\hk#1{\hskip #1mm}
\def\vk#1{\vskip #1mm}

\let\x=\times

\let\ph=\varphi

\let\wtld=\widetilde
\let\what=\widehat

\let\<=\langle
\let\>=\rangle
\let\lt=\left
\let\rt=\right

\let\otim=\otimes

\let\cv=\rightarrow

\let\disp=\displaystyle

\def\engl{\count99}
\def\chapno{\count100}
\def\parno{\count101}
\def\bibno{\count102}

\def\incr(#1,#2){\advance#1 by #2}

\long\def\bib#1#2#3{\item{[\the\bibno]}{\und{#1}}~{\sl #2}~{#3\par}
\incr(\bibno,1)}

\engl=0

\def\sect#1{{\bf\ifnum\chapno<1000{\the\chapno.}\fi\the\parno\ #1}
\global\incr(\parno,1)}

\def\lem#1{\ifodd\engl\sect{Lemma #1:}\else\sect{Lemme #1:}\fi\sl}
\def\thm#1{\ifodd\engl\sect{Theorem #1:}\else\sect{Th\'eor\`eme #1:}\fi\sl}
\def\cor#1{\ifodd\engl\sect{Corollary #1:}\else\sect{Corollaire #1:}\fi\sl}
\def\defi#1{\ifodd\engl\sect{Definition #1:}\else\sect{D\'efinition #1:}\fi\sl}
\def\defs#1{\ifodd\engl\sect{Definitions #1:}\else\sect{D\'efinitions #1:}\fi\sl}
\def\rem#1{\ifodd\engl\sect{Remark #1:}\else\sect{Remarque #1:}\fi}
\def\rems#1{\ifodd\engl\sect{Remarks #1:}\else\sect{Remarques #1:}\fi}
\def\prop#1{\sect{Proposition #1:}\sl}
\def\dm{\rm\ifodd\engl Proof: \else D\'emonstration: \fi}

\def\new{{\bf\ifnum\chapno<1000{\the\chapno.}\fi\the\parno\ }
\global\incr(\parno,1)}

\catcode`\$=3

\bibno=1
\chapno=1000
\parno=1

\let\nab=\nabla

\def\og{\leavevmode\raise.3ex\hbox{$\scriptscriptstyle\langle\!\langle$}}
\def\fg{\leavevmode\raise.3ex\hbox{$\scriptscriptstyle\,\rangle\!\rangle$}~}

\def\leqno{\eqno}

\engl=1
\hfill  \fdate{30}{avril}{2009}  \vk{1}

\hfil {\twbf REMARKS ON THE FRACTIONAL \par \hfil BROWNIAN MOTION}   \parag     \parag
\hfil {\bf Denis Feyel}  \vk{1}
D\'epartement de Maths., Universit\'e d'Evry-Val d'Essonne, Boulevard
Fran\c cois Mitterand, 91025 Evry cedex, France, {\tt denis.feyel@univ-evry.fr}  \vk{1}
\hfil {\bf Arnaud de La Pradelle}  \vk{1}
Laboratoire d'analyse Fonctionnelle, Universit\'e Paris VI, 4 place Jussieu,
75052 Paris, France, {\tt adelapradelle@free.fr}   \parag     \parag  

{\bf Abstract.} We study the fBm by use of convolution of the standard white
noise with a certain distribution. This brings some simplifications and new
results.   \parag     \parag  



   \parag     \parag  
{\bf Introduction}.   \parag  

There are many ways to tackle the fractional Brownian motion. In this paper,
we use a convolution of a white noise by a distribution $T$. This
distribution operates in principal value as explained in paragraphs I and
II. In paragraph III anf IV is defined a Skorohod type integral with respect
to the fBm. This allows to define vector valued rough paths which lead to
rough paths in the sense of T.Lyons. In paragraph V, we indicate a
regularization process of the fBm by convolution with some examples. We
study convergence of Riemann sums in paragraph VI, this also leads to
approximations by piecewise linear processes of fBm-Skorohod and
fBm-Stratonovich type integrals. Paragraph VII is devoted to prove that
every $\R^d$-valued fBm defined on $\R$ can be studied in this way.  \vk{1}
This work is intended to simplify many previous papers (loc.cit.), and
brings some new results.   \parag     \parag        \vfill  \eject  

{\bf I. The distribution $T$.}   \parag     \parag  
Consider the distribution-function
$$S(t)=\Frac{|t|^{\alpha -1}}{2\Gamma (\alpha )\cos(\pi \alpha /2)}=K(\alpha )|t|^{\alpha -1}$$
which is locally integrable for $\alpha >0, \alpha \dif{1}$. The distribution derivative
$$T(t)=\Frac{|t|^{\alpha -2}\op{Sign}(t)}{2\Gamma (\alpha -1)\cos(\pi \alpha /2)}$$
analytic with respect to $\alpha \in \;]1/2,3/2[$,
acts on the $H'$-H\"older continuous functions on $[a,b]$ with $\alpha +H'>1$
in the following way  \vk{1}
-- $\<T,\ph \>=\disp\int \ph (t)T(t)\,dt$ if $0\not\in [a,b]$,  \vk{1}
-- $\op{p.v.}\int _a^b\ph (t)T(t)\,dt=\ph (0)[S(b)-S(a)]+\disp\int _a^b[\ph (t)-\ph (0)]T(t)\,dt$
if $0\in ]a,b[$  \vk{1}
-- not defined if $0=a~\hb{or}~b$.   \parag  
Given $t$, the function $S_t(u)=S(u)-S(t-u)$ belongs to $L^2(du)$ for
$1/2<\alpha <3/2$. The Fourier transform is worth
$$\what S_t(\xi )=\int S_t(u)\ei{u\xi }\,du=\Frac{1-\ei{t\xi }}{|\xi |^\alpha }$$
so that
$$\<S_t,S_s\>_{L^2}=K(2\alpha )\lt[|t|^{2\alpha -1}+|s|^{2\alpha -1}-|t-s|^{2\alpha -1}\rt]$$
that is
$$\<S_t,S_s\>_{L^2}=K(2\alpha )\lt[|t|^{2H}+|s|^{2H}-|t-s|^{2H}\rt]\leqno  (1)$$
with $H=\alpha -1/2$ and $K(2\alpha )=[2\Gamma (2H+1)\sin(\pi H)]^{-1}$.   \parag  
Note that the map $\alpha \cv S_t(u)$ extends in a $L^2(du)$-valued holomorphic
function in the band $1/2<\op{Re}\alpha <3/2$. For $\alpha =1$, it is worth
$\pi ^{-1}\Log|1-t/u|$.    \parag  

\thm{} Let $\ph $ be a Banach-valued $H'$-H\"older continuous function defined
on a segment $[a,b]$, with $H+H'>1/2$. The convolution product
$$f(u)=\op{p.v.}\int _a^b\ph (t)T(u-t)\,dt\leqno  (2)$$
belongs to $L^2(du)$.  \vk{1}
\dm Note that $f$ is defined for $u\dif{a,b}$. First, for $u$ out of
$[a,b]$ we get
$$|f(u)|\le \|\ph \|_\infty \int _a^b|T(u-t)|\,dt=\|\ph \|_\infty |S(u-b)-S(u-a)|$$
Second, for $u\in \;]a,b[$, write
$\|\ph \|_{H'}=\Sup_{t\dif{s}}|\ph (t)-\ph (s)|/|t-s|^{H'}$, then
$$f(u)=\ph (u)[-S(u-b)+S(u-a)]+\int _a^b[\ph (t)-\ph (u)]T(u-t)\,dt$$
$$|f(u)|\le \|\ph \|_\infty |S(u-b)-S(u-a)|+K(\alpha ,H')\|\ph \|_{H'}\lt[|b-a|^{\alpha +H'-1}\rt]1_{]a,b[}(u)$$
with $K(\alpha ,H')=K(\alpha )/(\alpha +H'-1)$ (remark that $K(1,H')=2/(\pi H')$).
Finally we get
$$\eqalign{N_2(f)&\le \|\ph \|_\infty N_2(S_b-S_a)+K(\alpha ,H')\|\ph \|_{H'}|b-a|^{\alpha +H'-1/2}\cr
N_2(f)&\le [2K(2\alpha )]^{1/2}\|\ph \|_\infty |b-a|^H+K(\alpha ,H')\|\ph \|_{H'}|b-a|^{\alpha +H'-1/2}}$$
                                        \parag  
\rems{} 1{$^{\rm o}$}. If $\ph $ vanishes at a point of $[a,b]$, then
$\|\ph \|_\infty \le \|\ph \|_{H'}|b-a|^{H'}$, so that we get
$$N_2(f)\le K'(H,H')\|\ph \|_{H'}|b-a|^{H+H'}\leqno  (3)$$
2{$^{\rm o}$}. For $\alpha >1$, $T$ is locally integrable, so that the p.v. is unuseful.
   \parag     \parag  
{\bf II. The fBm}.   \parag     \parag  
Let $B_t$ be an $\R^d$ valued Brownian motion for $t\in \R$. As
$H=\alpha -1/2\in ]0,1[$, we have $S_t\in L^2(du)$ and we define
$$X^H_t=\int S_t(u)\,dB_u=K(\alpha )\int [|u|^{\alpha -1}-|t-u|^{\alpha -1}]\,dB_u\leqno  (4)$$
This is a centered Gaussian process with $H'$-H\"older continuous paths for
$H'<H$ (Kolmogorov lemma).     \vk{1}
The covariance matrix is given by formula (1)
$$\E(X^H_tX^H_s)_{ij}=\<S_t,S_s\>_{L^2}\delta _{ij}$$
Let $\ph $ be a $H'$-H\"older continuous function on a segment $[a,b]$, we define
$$\int _a^b\ph (t)\,dX^H_t=\int f(u)\,dB_u$$
where $f$ is the convolution defined in formula (2). Then formula (3) writes
$$N_2\lt(\int _a^b[\ph (t)-\ph (a)]\,dX^H_t\rt)\le K'(H,H')\|\ph \|_{H'}|b-a|^{H+H'}\leqno  (5)$$
The process $Y_t=\int _0^t\ph (t)\,dX^H_t$ centered Gaussian Banach valued process
for $H\in \;]0,1[$ and $H+H'>1/2$. Moreover the map $H\cv \disp\int _a^b\ph (t)\,dX^H_t$
is analytic in $H\in \;]1/2-H',1[$.   \parag     \parag  

{\bf Example:} Let $\wtld X^H$ an independant copy of $X^H$. If we take
$\ph (t)=\wtld X^H_t$ as a $L^2$-valued $H$-H\"older continuous function, we get
$$\int _a^b[\wtld X^H_t-\wtld X^H_a]\otim \,dX^H_t$$
By taking the coordinates, we can get the L\'evy areas of the fBm for $H>1/4$.
Putting for the $i$-coordinate
$$\int _a^b[X^{H,i}_t(\omega )-X^{H,i}_a(\omega )]\,dX^{H,i}_t(\omega )]=[X^{H,i}_b(\omega )-X^{H,i}_a(\omega )]^2/2$$
we get finally ageometric rough path calculus for $H>1/3$.    \parag     \parag  

{\bf III. Using the Skorohod integral}.   \parag     \parag  

Recall (\cf[4,5]) that the Gaussian Sobolev space $\cD^{1,2}(\Omega ,\mu )$
constructed upon the Gaussian measure $\mu $ is the space of Wiener functionals
$\Phi (\omega )$ such that $\nab{\Phi }(\omega ,\wtld \omega )$ belongs to $L^2(\mu \otim \mu )$. The divergence
is the transposed operator $\op{div}: L^2(\mu :\mu )\cv L^2(\mu )$.   \parag  
Let $\Phi $ be a $H'$-H\"older continuous $\cD^{1,2}$ valued function. Consider
the following integral
$$Z(\omega ,\wtld \omega )=\int _a^b\Phi _t(\omega )\otim \,dX^H_t(\wtld \omega )\leqno  (6)$$
then $Z$ belongs to $\cD^{1,2}\what {\otim }W^1$ where $W^1$ is the first Wiener chaos
and $\what {\otim }$ is the Hilbert tensor product. We get by formula (5)
$$\|Z-\Phi _a\otim (X^H_b-X^H_a)\|_{\cD^{1,2}\what {\otim }W^1}\le K(H,H')\|\Phi \|_{H',\cD^{1,2}}|b-a|^{H+H'}$$
Apply now the divergence operator, and define the fBm-Skorohod integral
\phantom{ab}
[1,6,7,8,13]
$$\int _a^b\Phi _t(\omega )\odot\,dX^H_t(\omega )=(\op{div}Z)(\omega )$$
Hence we get by continuity of the divergence
$$\lt\|\int _a^b[\Phi _t-\Phi _a]\odot\,dX^H_t\rt\|_{L^2(d\omega )}\le k'\|\Phi \|_{H,\cD^{1,2}}|b-a|^{H+H'}$$
   \parag  
{\bf Example:} suppose that $\Phi $ takes values into a non-homogneous Wiener
chaos of degree $k$ which is included in $\cD^{1,2}$ with an equivalent norm
as in $L^2$ [8]. For $H>1/4$ it is straightforward to verify that we can
define in this way
$$X^{H,(2)}_{ab,ij}=\int _a^bdX^H_{t,i}\odot\int _a^tdX^H_{s,j},\hk{6}\hb{and}\hk{6}
  X^{H,(3)}_{ab,ijk}=\int _a^bdX^H_{t,i}\odot{}X^{H,(2)}_{at,jk}$$
As $X^{H,(2)}_{ab}$ belongs to the second Wiener chaos, then $X^{H,(3)}_{ab}$
is well defined as an element of the third Wiener chaos for every $\{i,j,k\}$
(coordinate indices). Besides we have
$$\|X^{H,(2)}_{ab}\|_{L^2}\le \op{Cst}|b-a|^{2H},\hk{8}
  \|X^{H,(3)}_{ab}\|_{L^2}\le \op{Cst}|b-a|^{3H}$$
It is easily seen that we get a vector valued rough path for $H>1/4$,
for example we have for every $c\in [a,b]$
$$X^{H,(2)}_{ab}-X^{H,(2)}_{ac}-X^{H,(2)}_{cb}=(X^H_b-X^H_c)\otim (X^H_c-X^H_a)$$
Note the inversion of the tensor product, which does not matter.   \parag  
We could recover $\disp\int _a^bF(X^H_t)\odot\,dX^H_t$ thanks to the sewing
lemma ([11,12]). Observe that we can also obtain pathwise rough paths in the
sense of T.Lyons ([14]) for $H>1/4$.   \parag     \parag  

{\bf IV. Extending the classical calculus}.   \parag     \parag

Let $F$ be a polynomial. It is easily seen that $t\cv F(X^H_t)$ is H\"older
continuous, so that if $H>1/2$ we get
$$\int _a^bF(X^H_t)\odot\,dX^H_t=\op{div}\int _a^bF(X^H_t(\omega ))\otim dX^H_(\wtld \omega )$$
where the second member is a Young integral. By computing it, we find
$$\int _a^bF(X^H_t(\omega ))\,dX^H_t(\omega )-\int _a^b\nab{F}(X^H_t)\E(X^H_t\,dX^H_t)$$
and finally
$$\int _a^bF(X^H_t)\odot\,dX^H_t=\int _a^bF(X^H_t(\omega ))\,dX^H_t(\omega )-2HK(2\alpha )\int _a^b\nab{F}(X^H_t)t^{2H-1}\,dt$$
Then it is natural to put for $H>1/4$   \parag  
\defi{} For $H>1/4$ and $F$ a polynomial, we define
$$\int _a^bF(X^H_t(\omega ))\circ\,dX^H_t(\omega )=\int _a^bF(X^H_t)\odot\,dX^H_t+2HK(2\alpha )\int _a^b\nab{F}(X^H_t)t^{2H-1}\,dt$$
\rm
This formula can be read as an Ito formula for $H>1/4$.  Put

$$\wtld X^{H,(2)}_{ab,ij}=\int _a^bdX^H_{t,i}\rond \int _a^tdX^H_{s,j},\hk{6}\hb{and}\hk{6}
  \wtld X^{H,(3)}_{ab,ijk}=\int _a^bX^H_{tb,i}\rond \,dX^H_{t,j}\rond X^H_{at,k}$$
As for the Skorohod integral, we see that
$\wtld X^{H,(2)}_{ab}$ belongs to the second Wiener chaos, so that $\wtld X^{H,(3)}_{ab}$
is well defined as an element of the third Wiener chaos for every $\{i,j,k\}$.
Besides we have
$$\|\wtld X^{H,(2)}_{ab}\|_{L^2}\le \op{Cst}|b-a|^{2H},\hk{8}
  \|\wtld X^{(H,3)}_{ab}\|_{L^2}\le \op{Cst}|b-a|^{3H}$$
Now, for $H>1/2$ we get standard Young integrals, so that the rough path
algebraic relations are satisfied. They also hold for $H>1/4$ thanks to the
analyticity with respect to $H$. Hence, we get another vector valued
rough path for $H>1/4$.   \parag     \parag  

{\bf V. Approximations of the fBm}.   \parag     \parag  

{\bf 1{$^{\rm o}$}}. Let $\rho _n(t)\ge 0$ be a regularizing sequence, that is $\int \rho _n(t)\,dt=1$ and $\rho _n$
converges narrowly to the Dirac mass at 0. Put
$$S_n(u)=S*\rho _n(u),~~~~S_{t,n}(u)=S_t*\rho _n(u),~~~~~T_{t,n}(u)=T*\rho _n(u)$$
Put
$$X^H_{t,n}=\int S_{t,n}\,dB_u,~~~~~~~\int _a^b\ph (t)\,dX^H_{t,n}=\int f_n(u)\,dB_u$$
with
$$f_n(u)=\int _a^b\ph (t)T_{t,n}(u)\,dt$$
As $n$ goes to infinity, $S_{n,t}$ converges to $S_t$ in $L^2(du)$ so that
$X^H_{t,n}$ converge to $X^H_t$ in the first Wiener chaos.  \vk{1}
We have $f_n=\ph _{ab}*(T*\rho _n)$ where $\ph _{ab}$ is worth $\ph $ in $[a,b]$ and 0
elsewhere. It is easily seen that we have the associative relation
$f_n=(\ph _{ab}*T)*\rho _n$. As $\ph _{ab}*T$ belongs to $L^2$, the convolution by $\rho _n$ converges to $f$ in
$L^2$ when $n\cv \infty $. It follows that
$$\int _a^b\ph (t)\,dX^H_t=\Lim_{n\cv \infty }\int _a^b\ph (t)\,dX^H_{t,n}$$
   \parag  
{\bf Example:} take the 1-dimensional process $B_u$, and consider for $y>0$
$\rho _y(t)=\Frac{y}{\pi (t^2+y^2)}$. As well known $T_{t,y}(u)=\rho _y*T_t(u)$ is
harmonic with respect to $(t,y)\in \R\x\R^+$. Hence we get a harmonic extension
$X^H_{t,y}$ with respect to $(t,y)$. When $y$ goes to 0, $X^H_{t,y}$ converges
to $X^H_t$. Note that $\disp\int _a^b\ph (t)\,dX^H_{t,y}$ is a usual integral with
respect to a $\cC^\infty $-function $X^H_{t,y}$ for every $y>0$. It should be
observed that $X^H_{t,y}$ is the real part of a holomorphic function of
$t+\im{y}$ in the upper half-plane. This holomorphic function has been
investigated in [15].     \parag     \parag  

{\bf 2{$^{\rm o}$}}. (\cf [6,7]). For $\lambda >0$ put
$$G^\alpha _\lambda (t)=\Frac{2^{-\alpha /2}}{\Gamma (\alpha /2)}\int _0^\infty u^{\alpha /2-1}\e{-\lambda u}h_u(t)\,du$$
where
$$h_u(t)=(2\pi u)^{-1/2}\exp(-t^2/2u)$$
Then
$$\what {G^\alpha _\lambda }(\xi )=2^{-\alpha /2}(\lambda +\xi ^2/2)^{-\alpha /2}$$
Put
$$Y^{H,\lambda }_t=\int G^\alpha _\lambda (t-u)\,dB_u$$
$$\op{Cov}(Y^{H,\lambda }_t,Y^{H,\lambda }_s)=\int G^\alpha _\lambda (t-u)G^\alpha _\lambda (s-u)\,du=G^{2\alpha }_\lambda (t-s)=G^{2\alpha }_\lambda (h)$$
$$X^{H,\lambda }_t=Y^{H,\lambda }_t-Y^{H,\lambda }_0=\int [G^\alpha _\lambda (u)-G^\alpha _\lambda (t-u)]\,dB_u$$
For $\alpha -1/2\in ]0,1[$, we get by the Fourier transform
$$\op{Cov}(Y^{H,\lambda }_t,Y^{H,\lambda }_s)=\Frac{2^{-\alpha }}{2\pi }\int \Frac{\cos{h\xi }}{(\lambda +\xi ^2/2)^\alpha }\,d\xi $$
$$\E(X^{H,\lambda }_t,X^{H,\lambda }_s)=\Frac{2^{-\alpha }}{2\pi }\int \Frac{1-\cos{t\xi }-\cos(s\xi )+\cos{(t-s)\xi }}{(\lambda +\xi ^2/2)^\alpha }\,d\xi $$
As $\lambda \cv 0$, $G^\alpha _\lambda (u)-G^\alpha _\lambda (t-u)$ converges towards $S_t(u)$ in $L^2(du)$, so
that $X^{H,\lambda }$ converges towards $X^H$ in $L^2(d\omega )$.   \parag     \parag

{\bf VI. Riemann sums associated with a partition.}   \parag     \parag  

Let $\Delta =\{a=t_0,t_1,\ldots,t_n=b\}$ be a partition of $[a,b]$ with mesh $\delta $.
Choose points $\tau _i\in [t_i,t_{i+1}]$, and consider the Riemann sum
$$Z=\sum  _a^b\ph (\tau _i)[X^H_{t_{i+1}}-X^H_i]$$
where $\ph $ is $H'$-H\"older continuous on $[a,b]$.
Also consider the function
$$J(u)=\sum  _a^b\ph (\tau _i)[S_{t_{i+1}}(u)-S_i(u)]$$   \parag  

\lem{} For $H+H'>1/2$, $J(u)$ converges to $\ph _{ab}*T(u)$ in $L^2(du)$ as $\delta $
vanishes.   \vk{1}
\dm We have
$$\ph _{ab}*T(u)=\sum  _a^b\op{p.v.}\int _{t_i}^{t_{i+1}}\ph (t)T(u-t)\,dt$$
Hence the difference
$$D(u)=\ph _{ab}*T(u)-J(u)=\sum  _a^b\op{\;p.v.}\int _{t_i}^{t_{i+1}}[\ph (t)-\ph (\tau _i)]T(u-t)\,dt$$
For $u\not\in [a,b]$ we have
$$|D(u)|\le \|\ph \|_{H'}|S_b(u)-S_a(u)|\delta ^{H'}$$
For $u\in \,]a,b[$, we first look for $u\in ]t_k,t_{k+1}[$. We have
$$D(u)=\lt(\int _a^{t_k}+\int _{t_{k+1}}^b\rt)[\ph (t)-\ph (\tau _i)]T(u-t)\,dt\;+~
\op{p.v.}\int _{t_k}^{t_{k+1}}[\ph (t)-\ph (\tau _i)]T(u-t)\,dt$$
$$|D(u)|\le 2\|\ph \|_{H'}\delta ^{H'}[S_{t_k}(u)+S_{t_{k+1}}(u)]+|R(u)|$$
with
$$\eqalign{R(u)&=[\ph (u)-\ph (\tau _i)][S_{t_{k+1}}(u)-S_{t_k}(u)]+\int _{t_k}^{t_{k+1}}[\ph (t)-\ph (u)]T(u-t)\,dt\cr
|R(u)|&\le \|\ph \|_{H'}\delta ^{H'}[S_{t_k}(u)+S_{t_{k+1}}(u)]+\|\ph \|_{H'}\delta ^{\alpha +H'-1}/(\alpha +H'-1)\cr
|D(u)|&\le 3\|\ph \|_{H'}\delta ^{H'}[S_{t_k}(u)+S_{t_{k+1}}(u)]+\|\ph \|_{H'}\delta ^{\alpha +H'-1}/(\alpha +H'-1)}$$
  \vk{1} We get
$$\eqalign{\int D(u)^2\,du&=\lt(\int _{-\infty }^a+\int _b^\infty \rt)D(u)^2\,du+\sum  _a^b\int _{t_k}^{t_{k+1}}D(u)^2\,du\cr
\int D(u)^2\,du&\le K\|\ph \|_{H'}^2\delta ^{2H'}|b-a|^{2H}+K\|\ph \|_{H'}^2\sum  _a^b\delta ^{2H'+2\alpha -1}\cr
\int D(u)^2\,du&\le K\|\ph \|_{H'}^2[\delta ^{2H'}|b-a|^{2H}+\delta ^{2H+2H'-1}|b-a|]}$$
which converges to 0 as $\delta $ vanishes.   \parag  

\rem{} In fact this lemma also holds with a function $\tau _i(u)$ in place of
$\tau _i$.   \parag  
Then we can claim   \parag  

\thm{} We have in $L^2(d\omega )$
$$\int _a^b\ph (t)dX^H_t(\omega )=\Lim_{\delta \cv 0}\sum  _a^b\ph (\tau _i)[X^H_{t_{i+1}}-X^H_{t_i}](\omega )$$
\rm
{\bf Application to the fBm integrals.}   \parag  
Applying this result to equation (6), we get
$$Z(\omega ,\wtld \omega )=\Lim_{\delta \cv 0}\sum  _a^b\Phi _{\tau _i}(\omega )[X^H_{t_{i+1}}-X^H_{t_i}](\wtld \omega )$$
in the space $\cD^{1,2}(d\omega )\what {\otim }L^2(d\wtld \omega )$. By taking the divergence, we get
\hfuzz=2cm
$$\int _a^b\Phi _t(\omega )\odot\,dX^H_t(\omega )=\Lim_{\delta \cv 0}\sum  _a^b\Phi _{\tau _i}(\omega )[X^H_{t_{i+1}}-X^H_{t_i}](\omega )
-\wtld E[\nab{\Phi _{\tau _i}}(\omega ,\wtld \omega )[X^H_{t_{i+1}}(\wtld \omega )-X^H_{t_i}(\wtld \omega )]$$
in $L^2(d\omega )$. When $\Phi _t=F(X^H_t)$, we have
$$\nab{\Phi _t}(\omega ,\wtld \omega )=\nab{F}(X^H_t(\omega ))X^H_t(\wtld \omega )$$
$$\eqalign{&\sum  _a^b\nab{F}(X^H_{\tau _i}(\omega ))\wtld E[X^H_{\tau _i}(\wtld \omega )[X^H_{t_{i+1}}(\wtld \omega )-X^H_{t_i}(\wtld \omega )]=\cr
&=K(2\alpha )\sum  _a^b\nab{F}(X^H_{\tau _i}(\omega ))[|t_{i+1}|^{2H}-|t_i|^{2H}-|t_{i+1}-\tau _i|^{2H}+|t_i-\tau _i|^{2H}]\cr
&=2HK(2\alpha )\int _a^b\nab{F}(X^H_t(\omega ))t^{2H-1}\,dt-K(2\alpha )\sum  _a^b\nab{F}(X^H_{\tau _i}(\omega ))[|t_{i+1}-\tau _i|^{2H}-|t_i-\tau _i|^{2H}]}$$
\parag  
\thm{} Suppose that $\tau _i$ is the midpoint of $[t_i,t_{i+1}]$ then we have
$$\int _a^bF(X^H_t)\odot\,dX^H_t=\Lim_{\delta \cv 0}\sum  _a^bF(X^H_{\tau _i})[X^H_{t_{i+1}}-X^H_{t_i}]
-2HK(2\alpha )\int _a^b\nab{F}(X^H_t)t^{2H-1}\,dt$$
$$\int _a^bF(X^H_t)\rond \,dX^H_t=\Lim_{\delta \cv 0}\sum  _a^bF(X^H_{\tau _i})[X^H_{t_{i+1}}-X^H_{t_i}]$$
Hence the limite above exists in $L^2(d\omega )$.   \parag  \rm

{\bf Piecewise linear interpolation.}   \parag  
Let $\ph $ be a $H'$-H\"older continuous function, and
let $\ph _n$ be the function deduced from $\ph $ by linear interpolation
with vertices $t_i$ of $\Delta _n$. It is easily seen that the H\"older semi-norm
satisfies $\|\ph _n\|_{H'}\le \|\ph \|_{H'}$ and that $\|\ph -\ph _n\|_{H''}$ converges to 0 for
every $H''<H'$. Hence $J(\ph _n)$ is bounded in $L^2(du)$ and
$$N_2\lt(\int _a^b[\ph (t)-\ph _n(t)]\,dX^H_t\rt)\le k(H)\|\ph -\ph _n\|_{H''}|b-a|^{H+H'}$$
which tends to 0 when $n$ tends to infinity.   \parag  
Let $X^H_{t,n}$ be the linear interpolation of $X^H_t$ with respect to $\Delta _n$,
and consider $\ph _n=F(X^H_{t,n})$ for a polynomial $F$. As $\ph _n$ belongs to
a non-hogoneous Wiener chaos of degree $k$, we have for $s\le t$
$$\eqalign{F(X^H_{t,n})-F(X^H_{s,n})&=\int _0^1\nab{F}(X^H_{s,n}+\lambda [X^H_{t,n}-X^H_{s,n}])[X^H_{t,n}-X^H_{s,n}]\,d\lambda \cr
N_2(F(X^H_{t,n})-F(X^H_{s,n}))&\le \Sup_{\lambda \in \1}N_4(\nab{F}(X^H_{s,n}+\lambda [X^H_{t,n}-X^H_{s,n}]))N_4(X^H_{t,n}-X^H_{s,n})}$$
By the Nelson inequalities for non-homogeneous chaos (\cf [8]), we get
$$N_2(F(X^H_{t,n})-F(X^H_{s,n}))\le M'(a,b)N_2(X^H_{t,n}-X^H_{s,n})\le M(a,b)|t-s|^H$$
where $M'(a,b)$ and $M(a,b)$ depends only on $a$ and $b$.    \parag  

\lem{} Let $S_{t,n}$ be the interpolation of $S_t$ relative to $\Delta _n$. Then
for any $H'$-H\"older continuous $\ph $ with $H+H'>1/2$, we have
$$\ph _{ab}*T=\Lim_{\delta _n\cv 0}\int _a^b\ph _n(t)\,dS_{t,n}$$

\dm First, if $\ph $ is real valued, there exists functions
$u\cv \tau _i(u)\in [t_i,t_{i+1}]$ such that the second member is worth
$$\int _a^b\ph _n(t)\,dS_{t,n}=\sum  _a^b\ph (\tau _i(u))[S_{t_{i+1}}-S_{t_i}]$$
Applying lemma 4 and remark 5, we obtain
$$N_2\lt(\ph _{ab}*T-\int _a^b\ph _n(t)\,dS_{t,n}\rt)\le K\|\ph \|_{H'}\lt[\delta ^{2H'}|b-a|^{2H}+\delta ^{2H+2H'-1}|b-a|\rt]^{1/2}$$
Second, if $\ph $ is Banach valued, we get the same inequality thanks to the
Hahn-Banach theorem.   \parag  

\prop{} We have
$$\int _a^b\ph (t)\,dX^H_{t}(\omega )=\Lim_{\delta _n\cv 0}\int _a^b\ph _n(t)\,dX^H_{t,n}(\omega )$$
the limit being taken in $L^2(d\omega )$. Moreover if $F$ is a polynomial,
$$\int _a^bF(X^H_{t}(\omega ))\,dX^H_{t}(\wtld \omega )=\Lim_{\delta _n\cv 0}\int _a^bF(X^H_{t,n}(\omega ))\,dX^H_{t,n}(\wtld \omega )$$
in $L^2(d\omega \otim d\wtld \omega )$.       \vk{1}
\dm The first assertion follows from the lemma, and the second one is a
particular case in view of above.   \parag  

\cor{}
$$\int _a^bF(X^H_{t})\odot\,dX^H_{t}=\Lim_{\delta _n\cv 0}\int _a^bF(X^H_{t,n})\,dX^H_{t,n}-2HK(2\alpha )\int _a^b\nab{F}(X^H_{t,n})t^{2H-1}\,dt$$
\dm A straightforward computation.   \parag  \rm
\cor{}
$$\int _a^bF(X^H_{t})\circ\,dX^H_{t}=\Lim_{\delta _n\cv 0}\int _a^bF(X^H_{t,n})\,dX^H_{t,n}$$   \parag  
\rm
\rem{} In particular, the $\odot$-type and $\circ$-type rough paths are
obtained in this way (linear interpolation). Note that the $\circ$-type
rough path was obtained by ([2]). See also [13]. \parag   \parag

{\bf VII. Retrieving the initial Brownian motion}.   \parag     \parag  

Now, let $X^H_t$ be given a fBm, that is a continuous centered Gaussian
process defined on $\R$, with covariance
$$\E(X^H_sX^H_t)=K(2\alpha )[|t|^{2H}+|s|^{2H}-|t-s|^{2H}]$$
We suppose that $H\in ]0,1[$ and $\alpha =H+1/2$. Does there exists a standard
Brownian motion $B_t$ such that formula (4) holds ?   \parag  
First, let $A_t$ be the $\R^d$-valued standard Brownian motion defined on
$\R$, and let $Y_t$ be the fBm defined above from $A_t$
$$Y_t=\int [S^\alpha (u)-S^\alpha (u-t)]\,dA_u$$
Then
$$A_t=-\int [S^{2-\alpha }(u)-S^{2-\alpha }(u-t)]\,dY_u\leqno  (7)$$
Indeed, (7) makes sense as $u\cv |u|^{\alpha -2}[S^{2-\alpha }(u)-S^{2-\alpha }(u-t)]$ belongs
to $L^2(du)$ as easily seen by Fourier transform. Put
$\ph (v)=S_t^{2-\alpha }(v)$, we get
$$\int [S^{2-\alpha }(u)-S^{2-\alpha }(u-t)]\,dY_u=\int (\ph *T^\alpha )(v)\,dA_v=\int [T^\alpha *S^{2-\alpha }](v)\,dA_v$$
By Fourier transform, we check that $T^\alpha *S^{2-\alpha }](v)=-1_{[0,t](v)}$, so that we
are done.   \parag  
Let $\ph (t)$ be a $H$-H\"older continuous function with compact support, we have
seen that
$$\int _a^b\ph (u)\,dY_u=\Lim_{\delta \cv 0}\sum  _a^b\ph (\tau _i)[Y_{u_{i+1}}-Y_{u_i}]$$
for every partition of $[a,b]$ with mesh $\delta \cv 0$. It follows that
$$\int _a^b\ph (u)\,dX^H_u=\Lim_{\delta \cv 0}\sum  _a^b\ph (\tau _i)[X^H_{u_{i+1}}-X^H_{u_i}]$$
exists since the covariances are the same. Now take
$$\ph _R(u)=[S^{2-\alpha }(u)-S^{2-\alpha }(u-t)]1_{[-R,R]}(u)$$
As $R\cv +\infty $, we get
$$A_t=-\Lim_{R\cv \infty }\int _R^R\ph _R(u)\,dY_u$$
then
$$B_t=-\Lim_{R\cv \infty }\int _R^R\ph _R(u)\,dX^H_u$$
exists. Then we obtain a continuous centered Gaussian process $B_t$ on $\Omega $,
which is a Brownian motion, as easily checked by computing the covariance
function. Conversely one retrieves $X^H_t$ with the formula
$$X^H_t=\int [S^\alpha (u)-S^\alpha (u-t)]\,dB_u$$   \parag     \parag  

\hfil {\bf REFERENCES}   \parag  
\def\dlp{Feyel,D.; de La Pradelle,A.}

\bib{Alos,E.;Mazet;Nualart,D.}{Stochastic calculus with respect to fBm with
Hurst parameter lesser than $1/2$.}{Stoch.processes and applications, 86,
121-139 (2000).}  \vk{1}
\bib{Coutin,L.;Qian,Z.}{Stochastic Analysis, rough path Analysis, and
fractional Brownian motions.}{P.T.R.F., t.122, 108-140, (2002).}  \vk{1}
\bib{Decreusefond,L; \"Ust\"unel,A.S.}{Stochastic Analysis of the Fractional
Brownian Motion.}{Potential Analysis, {\bf 10}, 177-214, (1999).}  \vk{1}
\bib{\dlp}{Espaces de Sobolev gaussiens.}
{Ann.Inst.Fourier, t.39, fasc.4, 1989, p.875-908}    \vk{1}
\bib{\dlp}{Capacit\'es gaussiennes.}
{Ann. Inst. Fourier, t.41, f.1, p.49-76, 1991}         \vk{1}
\bib{\dlp}{Fractional integrals and Brownian processes.}{Publications de
l'Univ. Evry-Val-d'Essonne,(1996).}    \vk{1}
\bib{\dlp}{On fractional Brownian processes.}{Potential Analysis, {\bf 10},
273-288, (1999).}  \vk{1}
\bib{Feyel,D.}{Polyn\^omes harmoniques et in\'egalit\'es de Nelson.}{Publications de
l'Univ. Evry-Val-d'Essonne,(1995).}    \vk{1}
\bib{\dlp}{The FBM Ito formula through analytic continuation}
{Electronic Journal of Probability, Vol. 6, 1-22, Paper n{$^{\rm o}$}26 (2001).}  \vk{1}
\bib{\dlp}{Curvilinear integrals along rough paths.}{Preprint, (2003)}  \vk{1}
\bib{\dlp}{Curvilinear integrals along enriched paths}{Elec. J. Proba., 11, 34,
860-892, (2006).}    \vk{1}
\bib{Feyel,D.; de La Pradelle,A.; Mokobodzki,G.}{A non-commutative Sewing
Lemma.}{\par Electronic Communications in Proba., {\bf 13}, 24-35, (2008).}   \vk{1}
\bib{Hu.Y.}{Integral transformations and anticipative calculus for
fractional Brownian motions.} {Mem. Amer. Math. Soc. 175 (2005), no. 825,
viii+127 pp.}     \vk{1}
\bib{Lyons,T.J; Qian,Z.}{System Control and Rough Paths.}
{\par Oxford Science Publications, (2002).}     \vk{1}
\bib{Unterberger,J.}{Stochastic calculus for fractional Brownian motion with
Hurst parameter $H>1/4$; a rough path method by analytic extension.}
{\par Arxiv:maths/0703697 v1 [Math.PR] (March 2007).}

\bye